\documentclass[12pt, letterpaper, oneside]{article}

\usepackage[T2A]{fontenc}
\usepackage[utf8]{inputenc}
\usepackage{graphicx} 
\usepackage{setspace}
\usepackage{amsmath}
\usepackage{amssymb}
\usepackage{bbm}
\usepackage[left=30mm, top=30mm, right=30mm, bottom=30mm]{geometry}

\usepackage{hyperref}
\usepackage{color}

\newcommand{\RNumb}[1]{\uppercase\expandafter{\romannumeral #1\relax}}

\newcommand\smallO{
  \mathchoice
    {{\scriptstyle\mathcal{O}}}
    {{\scriptstyle\mathcal{O}}}
    {{\scriptscriptstyle\mathcal{O}}}
    {\scalebox{.7}{$\scriptscriptstyle\mathcal{O}$}}
  }
\newtheorem{lemma}{Lemma}
\newtheorem{theorem}{Theorem}

\title{On exponential recurrence of «Markov-up» processes}
\author{D.O. Kalikaeva\footnote{Moscow State University \&  Institute for Information Transmission Problems, Moscow, Russia; email: diana.kalikaeva@math.msu.ru} 
}
\date{}
\begin{document}
\maketitle

\begin{abstract}

The theory of «Markov-up» processes is being developed. This is a new class of stochastic processes with ``partial'' markovian features; it could also be called ``one-sided Markov''. Such a behavior may be found in the real world and in technical literature. So, naturally, phenomena of this type require special mathematical models. This paper follows the proposal of such a model published recently where the issue of recurrence and positive recurrence was investigated. In this article the exponential recurrence of the same process is established under appropriate assumptions, which potentially leads to an exponential rate of convergence to the invariant measure.
\medskip


\medskip


\end{abstract}

\section{Introduction}

The paper studies ``Markov-up'' processes. The idea of such a process  was first offered by Alexander Dmitrievich Solovyov in the end of the last century, but for a long time his idea was not implemented by anyone as a rigorous mathematical model. Recently one possible version of such a model was proposed in \cite{VV22}. In this model the process behaves as a Markov process during growth (stopping), and as soon as it begins to go down the character of its motion becomes more complex, namely, its behavior begins to depend on the entire trajectory of continuous decline. Note that the ``decision'' to turn downward is also markovian at the first moment. The formal description of the model will be presented in the next section. 

In general, what is interesting in any application is whether or not the process shows ergodic features. As any other process, such a model may be considered as markovian if the state space is properly extended; at the same time, this extension should be as ``small'' as possible, because otherwise the transition measures may become too singular with respect to each other. Another thing is that this new model admits a naturally embedded Markov chain; hence, the tools of semi-Markov processes may be applied here. 

Some informal practical examples may be found in the literature and in the real world where some devices or networks show a behavior which looks markovian in one direction and much more complicated in the other direction. One example is the article \cite{IE} in which the authors studied  blackouts in North America and ways to improve the reliability of the power system. The power grid is an interconnected system and one failure can lead to the failure of all components, this phenomenon is called the cascade effect. So, the cause of a major blackout in North America in 1996 was cascading failures. The failure of two lines led to the tripping of the generator at McNary Dam, which led to the shutdown of more than a dozen power units. As a result of the chain reaction, 11 US states and 2 Canadian provinces were de-energized. 

Another example is the collapses in Global Stock Markets. The authors of the article \cite{M9} modeled local, regional and global crashes and came to the conclusion that global crashes do not start suddenly, they occur according to the domino principle: everything starts with local markets, then they infect regional ones. That is, minor crashes today significantly increase the probability of a major crash tomorrow.

In \cite{VV22} the authors consider the random variables $\tau$, characterizing the first moment when the process crosses some floor $N \in \mathbb Z_+$ and $\gamma$ - the shortest time when the process leaving $[0,N]$ took the value $N$. Under some assumptions, the results on the boundedness of $\mathbb{E}_{x}\tau$ and $\mathbb{E}_{x} \gamma$ are obtained. The results are used to prove the existence of an invariant measure for the process.

The purpose of this paper is to find assumptions that will allow to establish the finiteness of the exponential moment $\mathbb{E}_{x}e^{\alpha\tau}$ for some positive $\alpha$. The result will further allow to estimate the rate of convergence to the invariant measure.

\section{The model and the assumptions}

Let us consider the process $ X_n,$ $ n\geq{0}$ on $\mathbb{Z_+}$ (or on ${\mathbb{Z}_{0,\widetilde{N}}} = \{0, \dots, \widetilde{N}\}$ for some $\widetilde{N} \in \mathbb{N}$, $ \widetilde{N}<{\infty}$) and define random variables for $N \geq{0}$

\begin{equation}\label{zeta}
\zeta_n: = \inf(k\leq{n}:\Delta X_i = X_{i+1} - X_i<{0},  \forall i = k, \dots, n-1),
\end{equation}
\begin{equation}\label{xi}
\xi_n: = \sup(k\geq{n}:\Delta X_i = X_{i+1} - X_{i} \geq{0},  \forall i = n, \dots, k-1)\vee n,
\end{equation}
\begin{equation}\label{chi}
\chi_n: = \sup(k\geq{n}:\Delta X_i = X_{i+1} - X_{i}<{0}, \forall i = n, \dots, k-1)\vee n,
\end{equation}
\begin{equation}\label{tau}
\tau := \inf(t\geq{0}:X_t\leq{N}).
\end{equation}
Also, let
\begin{align}\label{tf}
   \hat{X}_{i,n} := X_i \mathbbm{1} (\min(\zeta_n, n)\leq{i}\leq{n}),\quad\widetilde{\mathcal{F}_n} := \sigma(\zeta_n; \hat{X}_{i,n}: 0\leq{i}\leq{n}).
\end{align}

Note that $\widetilde{\mathcal{F}_n}$ is not a filtration. And for all $n \in \mathbb{Z_+} $ $\widetilde{\mathcal{F}_n}\subset \mathcal{F}_n$ where $(\mathcal{F}_n)_{n \in \mathbb{Z_+}}$ is the natural filtration associated to the process $(X_n)_{n \in \mathbb{Z_+}}$.  

As usual, all inequalities with conditional probabilities and conditional expectations are understood a.s.. The following assumptions are made:

\textbf{A1. Random memory depth:}
\textit {For any $ n \geq{0} $}
        
\begin{equation}\label{ctF}
\mathsf P(X_{n+1} = j | \mathcal{F}_n) = \mathsf P(X_{n+1} = j | \widetilde{\mathcal{F}_n}) ,
\end{equation}
This condition distinguishes the proposed process from Markov chain. If $n$ is such that $X_{n-1} \le X_{n}$, then $\widetilde{\mathcal{F}}_{n} = \sigma (X_{n})$ and it may be informally said that at this moment the process has a  Markov property. Yet, if $X_{n-1} > X_{n}$, then $\widetilde{\mathcal{F}}_{n} = \sigma (\zeta_{n}, X_{\zeta _{n}}, \dots, X_{n} )$ and the conditional distribution of the next value $X_{n+1}$ depends on some nontrivial part of the past; informally it may be said that at this moment the behavior of the process is not markovian.
    
\textbf{A2. Irreducibility (local mixing):}
\textit { For any $x \leq{N} $ and for all two states $ y = x, y = x+1$}
\begin{center}
$\mathsf P(X_{n+1} = y | \widetilde{\mathcal{F}_n}, X_n = x) \geq{\rho} >{0}.$
\end{center}
Note that $2\rho \le 1$. The assumption A2 will guarantee the irreducibility of the process in the extended state space where the process becomes Markov.
        
\textbf{A3. Recurrence-1:}
\textit {There exists $N \geq{0}$ such that for the conditional probability of a jump down the following is performed}
\begin{equation}\label{q0}
\mathsf P(X_{n+1}<{X_n}|\widetilde{\mathcal{F}_n}, N<{X_n})\geq{\kappa_0}>{0},
\end{equation}
$$\mathsf P(X_{n+1}<{X_n}|\widetilde{\mathcal{F}_n}, N<{X_n}<{X_{n-1}})\geq{\kappa_1}>{0},$$ 
\textit{etc., for all $ n\geq{m} $}
\begin{equation}\label{q}
\mathsf P(X_{n+1}<{X_n}|\widetilde{\mathcal{F}_n}, N<{X_n}<{\dots}<{X_{n-m+1}})\geq{\kappa_{m-1}}>{0},\forall m\geq{1}.
\end{equation}
Note that $\{\kappa_n\}_{n = 0} ^{\infty}$ is a non-decreasing sequence.

Denote $q = 1 - \kappa_0 < 1$. Then for the conditional probability of a jump up we have
\begin{center}
$ \mathsf P(X_{n+1}\ge{X_n}|\widetilde{\mathcal{F}_n}, N<{X_n})\le 1 - {\kappa_0} = q < 1$.
\end{center}
    
\textbf{A4. Recurrence-2:}
\textit { It is assumed that the sequence $\{\kappa_i\}_{i=0}^\infty$ tends to 1 with exponential rate as $n \to \infty$.}

\textit {It is assumed that the following infinite product converges}
\begin{align}\label{prod}
\bar\kappa_{\infty} := \prod\limits_{i=0}^{\infty}\kappa_i > 0.
\end{align}
This assumption implies that $\kappa_i \to 1$ as $i\to\infty$. Given that $\{\kappa_i\}_{i = 0} ^{\infty}$ is monotonic, tends of sequence to 1 may be interpreted as follows: the longer the process is falling down, the higher the probability that it will continue to fall.

Let $\bar q = 1 - \bar\kappa_{\infty} < 1$. This is the upper bound for probability that in one
go the process will not reach $[0, N]$.
        
\textbf{A5. Exponential moment of the value of jump up is limited:}
\textit {For sufficiently small $\alpha>{0}$}
\begin{align}\label{emoment} 
M_{\alpha}:=\mathop{\mbox{ess}}\limits_{\mathbb P}\sup\limits_{\omega} \sup\limits_{n} \mathsf{E}(e^{\alpha(X_{n+1}-X_n)_+}|\widetilde{\mathcal{F}_n})<{\infty}.
\end{align}

\section{The existence of processes and an example}

Formally, our Markov-up-to-be process may be presented as the last component of the vector valued Markov process $(V_n)$ in the state space of a variable dimension $V_n\in (\mathbb Z_+)^{n+1}$, 
\begin{align*}
V_0=X_0, \; \ldots, \; V_n = (X_0, X_{1}, \ldots, X_{n-1}, X_{n}). 
\end{align*}
Denote 
$$
v_n = (x_0, \ldots, x_{n}), \quad x_i\in \mathbb Z_+, \quad 0\le i\le n.
$$
It is assumed that the evolution of $V_n$ at time $n$ obeys a probability transition kernel $$Q_n(v_n, v_{n+1}) = \mathsf P(V_{n+1} = v_{n+1} \vert V_n = v_n)$$ satisfying certain restrictions. Note that whatever these restrictions are, the measure on the space of trajectories of $(V_n)$ is well-defined due to Kolmogorov's theorem about consistent family of distributions:
$$
\mu_n(v_n) := \mathsf P_{v_0}(V_n = v_n) 
= \prod_{i=1}^{n}Q_i(v_{i-1}, v_i).
$$
in other words, with any set of transition kernels $(Q_n)$, the process $V_n$ exists. Recall that $X_n$ is defined as the first component of the vector $V_n$; hence, it also does exist. The process $X$ may be called a Markov chain with a variable memory.

What makes this process a Markov-up one is the following set of restrictions on the transition kernels $Q_n$. 
\begin{enumerate}
\item[(B1)]\label{it1}
The value $Q_n(v_n, v_{n+1})$ does not depend on the tail $(x_0, \ldots, x_{n-1})$ for any such vector with $Q_{n-1}(x_0, \ldots, x_{n-1})>0$ on the set $x_{n+1} \ge x_n$. In other words, 
$$
1(X_{n+1} \ge X_n) \mathsf P(X_{n+1} \in \cdot \, \vert V_n) =
1(X_{n+1} \ge X_n) \mathsf P(X_{n+1} \in \cdot \, \vert X_n).
$$

\item[(B2)]\label{it2}
For any vector $v_n = (x_0, \ldots, x_n)$ denote 
$$
z_n = \min(0\le k < n: x_{n}<x_{n-1}, \ldots , x_{n-k} < x_{n-(k+1)}); \quad \min(\emptyset) = -1.
$$
It is assumed that on the set $(z_n=k)$ for $k< n$ the value $Q_n(v_n, v_{n+1})$ does not depend on the tail $(x_0, \ldots, x_{n-(k+1)})$ for any such vector with $Q_{k-1}(x_0, \ldots, x_{k-1})>0$. In other words, 
$$
1(\hat\zeta_{n} = k) \mathsf P(X_{n+1} \in \cdot \, \vert V_n) =
1(\hat\zeta_{n} = k) \mathsf P(X_{n+1} \in \cdot \, \vert X_n, \ldots, X_{n-(k+1)}),  
$$
where
$$
\hat\zeta_n = \min(0 \le k < n: X_{n}<X_{n-1}, \ldots , X_{n-k}<X_{n-(k+1)}); \quad \min(\emptyset) = -1.
$$
\end{enumerate}

So, for such a model the process $X_n$ is regarded as the last component of the vector-valued Markov chain $V_n$ in the state space of a variable dimension, 
$$
V_n = (X_{n-(\hat\zeta_n+1)}, \ldots X_n).
$$
Equivalently  but, perhaps, a little more explicitly the process $(X_n)$ may be regarded as the component of the vector-valued Markov process $(V_n)$ accompanied by one more coordinate $\hat\zeta_n$ (although, the knowledge of the dimension of $V_n$ already provides this information): 
$$
\bar V_n = (\hat\zeta_n, X_{n-(\hat\zeta_n+1)}, \ldots, X_n).
$$
Both versions are also equivalent to the presentation of $X_n$ via the formula (15) in \cite{VV22}. Either representation clearly shows the existence of the process $(X_n)$.

\medskip

Here is an {\bf example} of a process satisfying the conditions (B1) -- (B2); in what follows we show that it also satisfies the assumptions (A1) -- (A5). It is a perturbation of a basic homogeneous Markov chain with transition probabilities $p_x=\mathsf P(X_1=x+1\vert X_0=x)=1/5, \, q_x=\mathsf P(X_1=x-1\vert X_0=x)=3/5, \, r_x = \mathsf P(X_1=x\vert X_0=x)=1/5$ for any $0<x\in \mathbb Z_+$, except the state $x=0$ where $p_0=\mathsf P(X_1=1\vert X_0=0)=1/2, \, r_0 = \mathsf P(X_1=0\vert X_0=0)=1/2$. Note that this process 
 
is ergodic, irreducible, and positive and exponentially recurrent with a unique invariant measure. Let us perturb the transition probabilities -- except only at the zero state -- as follows. The perturbed transition probabilities now depend on some ``past'': 
\begin{align}\label{jdown} 
&\mathsf P(X_{n+1}=x_n-1\vert V_n=v_n; \zeta_n=k, X_n=x_n) = \frac35 +  \frac25(\frac12 + \ldots + \frac1{2^k}), 
\end{align}
\begin{align}\label{stay} 
&\mathsf P(X_{n+1}=x_n\vert V_n=v_n; \zeta_n=k, X_n=x_n)  = \frac15 - \frac1{5}(\frac12 + \ldots + \frac1{2^k}) = \frac1{5\times2^k}, 
\end{align}
\begin{align}\label{jup} 
&\mathsf P(X_{n+1}=x_n+1\vert V_n=v_n; \zeta_n=k, X_n=x_n) =  \frac15 - \frac1{5}(\frac12 + \ldots + \frac1{2^k}) = \frac1{5\times2^k}.
\end{align}
The probability measure with such transition probabilities does exist and it satisfies (B1)--(B2) (in particular, the process satisfies (A1)) and  all other assumptions of the paper with the value $N=0$. Note, in particular, that 
\begin{align}\label{est} 
&\frac35\le \mathsf P(X_{n+1}\!=\!x_n-1\vert V_n=v_n; \zeta_n=k, X_n=x_n) \!<\! 1, 
 \nonumber \\\\ \nonumber 
&1\!-\!\mathsf P(x_{n+1}\!=\!x_n\!-\!1\vert v_n; z_n\!=\!k, X_n=x_n)\! =\! \frac2{5\times2^k}.
\end{align}

Condition (A2) is clearly verified because by the construction 
$$
\mathsf P(X_{n+1}=1\vert \widetilde {\cal F}_n, X_n = 0) 
= \mathsf P(X_{n+1}=0\vert \widetilde {\cal F}_n, X_n = 0) = \frac12.
$$

Let us show that the assumption (A3) is also true for the process. We have,

\begin{align*}
&\mathsf P(X_{n+1}<{X_n}|\widetilde{\mathcal{F}_n}, N<{X_n}<{\dots}<{X_{n-(m+1)}}) 
 \\\\
&= \mathsf{E}(\mathbbm{1}(X_{n+1}<X_n)|{{V}_n}, \hat{\zeta}_n = m , N<{X_n}<{\dots}<{X_{n-(m+1)}}) 
 \\\\
&=\mathsf{E}(\mathsf{E}(\mathbbm{1}(X_{n+1}<X_n)|{{V}_n}, \hat\zeta_n = m )|{{V}_n}, \hat\zeta_n = m , N<{X_n}<{\dots}<{X_{n-(m+1)}})
 \\\\
&=\mathsf{E}(\mathsf P(X_{n+1}=X_n-1|{{V}_n}, \hat\zeta_n = m )|{{V}_n}, \hat\zeta_n = m , N<{X_n}<{\dots}<{X_{n-(m+1)}})
 \\\\
&=\frac35 +  \frac25(\frac12 + \ldots + \frac1{2^m}).
\end{align*}
So, the justification of (A3) follows from the inequality
\begin{align*}
\kappa_{m+1} = \frac35 +  \frac25(\frac12 + \ldots + \frac1{2^m}) = 1 - \! \frac2{5\times2^m} > 0.
\end{align*}

Note that the sequence $\{\kappa_n\}_{n = 0}^\infty$ is increasing and tends to 1. Let us show that 
assumption (A4) holds for some $\alpha$. 
Consider the series
\begin{align*}
&\sum\limits_{k=1}^{\infty}e^{\alpha k}(1-\kappa_k) = \sum\limits_{k=1}^{\infty}e^{\alpha k} \frac2{5\times2^k} =  \frac25 \sum\limits_{k=1}^{\infty} \frac{e^{\alpha k}}{2^k}.
\end{align*}
Due to Cauchy's criterion for convergence we have that if
$e^{\alpha}/2 < 1$,  i.e. if $\alpha < \ln{2}$, then the series converges. The convergence of an infinite product $\displaystyle  \prod\limits_{k=0}^{\infty}\kappa_k = \prod\limits_{k=0}^{\infty} (1- \frac2{5\times2^k} )$ follows from the convergence of the series $\displaystyle\sum\limits_{k=0}^{\infty} \ln{(1- \frac2{5\times2^k} )}$ ( in this case we have that $\displaystyle
\ln{(1- \frac2{5\times2^k} )} = -\frac{2}{5\times2^k} + \smallO (\frac{1}{2^{2k}}) \text{ and } \frac12 < 1$). Thus, indeed, (A4) holds.

Also for $\alpha < \ln{2}$,
\begin{align*}
M_{\alpha}:=\mathop{\mbox{ess}}\limits_{\mathbb P}\sup\limits_{\omega} \sup\limits_{n} \mathsf{E}(e^{\displaystyle \alpha(X_{n+1}-X_n)_+}|V_n)< 1 + e^{\displaystyle\ln{2}} = 3 < \infty.
\end{align*}
Hence, the process satisfies the assumption (A5). 

\section{Auxiliary lemmata}
Below the standard notations $\mathsf P_{v_0}(\cdot) = \mathsf P(\cdot\vert V_0=v_0)$ and similarly for the expectation, $\mathsf E_{v_0}$, from the theory of Markov processes will be used.
\begin{lemma}\label{lemma1}
Under the assumption (A3) for any $ v_0>{N}$ and $\alpha < \ln{\frac{1}{q}},$ 
\begin{equation}
\mathsf{E}_{v_0}e^{\alpha(\xi_n-n)}\leq{ M_1^{(\alpha)}:={e^\alpha q\over{1-e^\alpha q}}}.		
\end{equation}
\end{lemma}
\textbf{Proof.}
For all ${i \geq{n}}$ let
\begin{center}
$e_i = \mathbbm{1}(X_{i+1}\geq{X_i})$,  $\bar{e}_i = \mathbbm{1}(X_{i+1}\textless{X_i})$,  $l^i_n = \bar{e}_i\prod\limits_{k=n}^{n-1}e_k$ $\Delta{X_i} = X_{i+1} - X_i$.
\end{center}
We have, for all ${i \geq{n}}$
\begin{center}
$
\mathsf{E}_{v_0}(e_i|X_i\textgreater{N}) = \mathsf{P}_{v_0}(X_{i+1}\geq{X_i}|X_i\textgreater{N}) $ 
$$
= \mathsf{E}_{v_0}(\mathbb{P}_{x}(X_{i+1}\geq{X_i}| \widetilde{\mathcal{F}_i}, X_i\textgreater{N})| X_i\textgreater{N}) \leq{1-\kappa_0} = q.
$$
\end{center}
Then almost surely
\begin{center}
$e^{\alpha(\xi_n-n)} =  \sum\limits_{k=0}^{\infty}e^{\alpha k}\bar {e}_{n+k}\prod\limits_{i=n}^{n+k-1}e_i = \sum\limits_{k=1}^{\infty}e^{\alpha k}\bar {e}_{n+k}\prod\limits_{i=n}^{n+k-1}e_i = \sum\limits_{k=1}^{\infty}e^{\alpha k}l^{n+k}_n.$
\end{center}
Using this representation, get
\begin{center}
$\mathsf{E}_{v_0}e^{\alpha(\xi_n-n)} = \mathsf{E}_{v_0}\sum\limits_{k=1}^{\infty}e^{\alpha k}\bar {e}_{n+k}\prod\limits_{i=n}^{n+k-1}e_i \leq{}\sum\limits_{k=1}^{\infty}e^{\alpha k}\mathsf{E}_{v_0}\prod\limits_{i=n}^{n+k-1}e_i\leq{}\sum\limits_{k=1}^{\infty}e^{\alpha k}q^k$
                
$$ =\sum\limits_{k=1}^{\infty}(e^{\alpha}q)^k = {e^\alpha q\over{1-e^\alpha q}} =: M_1^{(\alpha)} \textless{\infty}.$$
\end{center}
Lemma 1 is proved. \hfill $\square$

\begin{lemma}\label{lemma2}
Under the assumptions (A3) and (A4) $\forall$ $ v_0\textgreater{N},(n=0)$ 
\begin{equation}
\mathsf{E}_{v_0}e^{\alpha(\chi_n-n)}\mathbbm{1}(\chi_n\textless{\tau})\leq\sum\limits_{i=1}^{\infty}e^{\alpha i}(1-\kappa_i)=:M_2^{(\alpha)}\textless{\infty}.
\end{equation}
\end{lemma}
        \textbf{Proof.}
        As in the previous lemma , in the same notation we have
            \begin{center}
                 $e^{\alpha(\chi_n-n)}\mathbbm{1}(\chi_n\textless{\tau}) \leq{}  \sum\limits_{k=1}^{\infty}e^{\alpha k}e_{n+k}\mathbbm{1}(n+k-1\textless{\tau})\prod\limits_{i=n}^{n+k-1}\bar{e}_i$.
            \end{center}
        So, 
              $$\mathsf{E}_{v_0}e^{\alpha(\chi_n-n)}\mathbbm{1}(\chi_n\textless{\tau}) 
              \leq{}  \mathsf{E}_{v_0}\sum\limits_{k=1}^{\infty}e^{\alpha k}e_{n+k}\mathbbm{1}(n+k-1\textless{\tau})\prod\limits_{i=n}^{n+k-1}\bar{e}_i$$
            $$\leq{}\sum\limits_{k=1}^{\infty}e^{\alpha k}\mathsf{E}_{v_0}\mathbbm{1}(n+k-1\textless{\tau})(\prod\limits_{i=n}^{n+k-1}\bar{e}_i) \mathsf{E}_{v_0}(e_{n+k}|\Delta{X_i}\textless{0}, n\leq{}i\leq{}n+k-1) $$
              $$\stackrel{A3}{\leq{}}\sum\limits_{k=1}^{\infty}e^{\alpha k}\mathsf{E}_{v_0}\mathbbm{1}(n+k-1\textless{\tau})(\prod\limits_{i=n}^{n+k-1}\bar{e}_i)(1-\kappa_k) \leq{}\sum\limits_{k=1}^{\infty}e^{\alpha k}(1-\kappa_k) =:M_2^{(\alpha)} \stackrel{A4}{\textless{}} \infty.$$
         Lemma 2 is proved. \hfill $\square$ 

    \begin{lemma}\label{lemma3}
		Under the assumptions  (A3) and (А5) for all $ v_0\textgreater{N},$ 
		\begin{equation}
         \sup_{n,x}\mathsf{E}_{v_0}e^{\alpha (X_{\xi_n}-X_n)}\leq{M_3^{(\alpha)}}\textless{\infty}.
		\end{equation}
	\end{lemma}
        \textbf{Proof.}
        In the same notation we have
        \begin{center}
            $\mathsf{E}_{v_0}e^{\alpha (X_{\xi_n}-X_n)} = \mathsf{E}_{v_0}\sum\limits_{i=n+1}^{\infty} l^i_n e^{\alpha(X_i-X_n)}\leq{} \sum\limits_{i=n+1}^{\infty} \mathsf{E}_{v_0}(\prod\limits_{j = n}^{i-1}e_j)e^{\alpha \sum\limits_{k=n}^{i-1}\Delta{X_k} }=\sum\limits_{i=n+1}^{\infty}\mathsf{E}_{v_0} \prod\limits_{j = n}^{i-1}(e_j e^{\alpha\Delta X_j}).$       
        \end{center}
        Lets consider general term of the series. Applying Cauchy-Bunyakovsky-Schwarz inequality, get estimation
        \begin{center}
            $ \mathsf{E}_{v_0} (\prod\limits_{j = n}^{i-1}e_j e^{\alpha\Delta X_j}) = \mathsf{E}_{v_0}\mathsf{E}_{\mathcal{F}_{i-1}} (\prod\limits_{j = n}^{i-1}e_j e^{\alpha\Delta X_{j}}) = \mathsf{E}_{v_0}  (\prod\limits_{j = n}^{i-2}e_j e^{\alpha\Delta X_{j}})\mathsf{E}_{\mathcal{F}_{i-1}}e_{i-1} e^{\alpha\Delta X_{i-1}}$

            $  \stackrel{K-B} {\leq{}}\mathsf{E}_{v_0}  (\prod\limits_{j = n}^{i-2}e_j e^{\alpha\Delta X_{j}})(\mathsf{E}_{\mathcal{F}_{i-1}}e_{i-1})^{1\over{2} }(\mathsf{E}_{\mathcal{F}_{i-1}}e^{2\alpha\Delta X_{i-1}})^{1\over{2}}$
            
            $\stackrel{A5} {\leq{}} q^{1\over{2}} M^ {1\over{2}}_{2\alpha}\mathsf{E}_{v_0}  (\prod\limits_{j = n}^{i-2}e_j e^{\alpha\Delta X_{j}})\leq{}(q^{1\over{2}} M^ {1\over{2}}_{2\alpha})^{i-n}.$
            
        \end{center}
        The last inequality is obtained inductively.
        Let $  \mu_{\alpha}:= q^{1\over{2}} M^ {1\over{2}}_{2\alpha}$, for sufficiently small $\alpha$ we have $\mu_{\alpha} \textless{1}$.
        So,
        \begin{center}
            $\mathsf{E}_{v_0}e^{\alpha (X_{\xi_n}-X_n)}\leq{}\sum\limits_{i=n+1}^{\infty} \mu_{\alpha}^{i-n} = \sum\limits_{i=1}^{\infty} \mu_{\alpha}^{i} = {\mu_{\alpha}\over{1- \mu_{\alpha}}} =: M_3^{(\alpha)}\textless{\infty}.$
        \end{center}
         Lemma 3 is proved. \hfill $\square$ 
\section{Main results}
       Let us recall,

        \begin{center}
        $\tau := \inf(t\geq{0}:X_t\leq{N})$.
        \end{center}
        
    \begin{theorem}\label{teorema}
		Under the assumptions (A1), (A3)-(А5)  for sufficiently small $\alpha \textgreater{0}$ there exist constants $C_1\textgreater{0}$ such that 
        \begin{equation}
        \mathsf{E}_{v_0}e^{\alpha \tau} \leq{e^{\alpha v_0} C_1.}
		\end{equation}
	\end{theorem}
        \textbf{Proof.}
        If $v_0\leq{N} $, then $\tau = 0$, this case is trivial, therefore, we assume that $v_0\textgreater{N}$.
        Consider the events:
        \begin{center}
            $ A_i = \{$exactly $i-1$ unsuccessful attempts to descend to the floor $[0, N],$ attempt №$ i$ is successful $\}, i\geq{1},$

            $ B_j = \{$attempt №$ j$ to reach $[0, N]$ is unsuccessful$\}, j\geq{1},$
            
            $ B_j^c = \{$attempt №$ j$ to reach $[0, N]$ is successful$\}, j\geq{1}.$
        \end{center}
        Note that ( according to the assumption А3) the probability of an unsuccessful attempt to cross the floor (that is , event $ B_j$) less then $\bar q$. In this notations for event $A_i$ is valid: $\tau = T_i$, $A_i = (\cap_{j=1}^{i-1} B_j)\cap B_i^c$, $P(A_i) \leq{\bar {q}^{i-1}}$.

\medskip

        \textbf{Case \RNumb{1}}: at $t=0$ the process is falling down.
        
        Let us define stopping times:
        \begin{center}
            $t_0 = T_0 = 0,$ $ T_1 = \chi_{t_0},$ $t_1 = \xi_{T_1},$ $ T_2 = \chi_{t_1},$ $t_2 = \xi_{T_2},$ $ T_3 = \chi_{t_2}, \dots$
        \end{center}
        So, from $t_{i-1}$ to $ T_i$ the process is continuously falling, at $T_i$ the fall is replaced by growth and up to $ t_i$ the process continuously run up. Process a.s. finite number of times will change its behavior until it reaches $[0, N]$.
        Note that $T_i-t_{i-1}\leq{X_{t_{i-1}}}$ и $B_j \in \mathcal F_{T_j}$.
        
        Let estimate
        \begin{center}
            $\mathsf{E}_{v_0}e^{\alpha \tau} =
            \sum\limits_{i \geq{} 1} \mathsf{E}_{v_0}e^{\alpha \tau}\mathbbm{1}(A_i)= 
            \sum\limits_{i \geq{} 1} \mathsf{E}_{v_0}e^{\alpha T_i}\mathbbm{1}((\cap_{j=1}^{i-1}B_j)\cap B_i^c)\stackrel{2}{=} 
            \sum\limits_{i \geq{} 1} \mathsf{E}_{v_0}(\prod\limits_{j=1}^{i-1}\mathbbm{1}(B_j))\mathbbm{1}(B_i^c)e^{\alpha T_i}$
            
            $\stackrel{3}{=}
            \sum\limits_{i \geq{} 1}\mathsf{E}_{v_0}\mathsf{E}_{\mathcal{F}_{t_{i-1}}}(\prod\limits_{j = 1}^{i-1}\mathbbm{1}(B_j))\mathbbm{1}(B_i^c)e^{\alpha T_i}\stackrel{4}{=}
            \sum\limits_{i \geq{} 1}\mathsf{E}_{v_0}(\prod\limits_{j = 1}^{i-1}\mathbbm{1}(B_j))\mathsf{E}_{\mathcal{F}_{t_{i-1}}}\mathbbm{1}(B_i^c)e^{\alpha T_i}$

            $\stackrel{5}{=}\sum\limits_{i \geq{} 1}\mathsf{E}_{v_0}(\prod\limits_{j = 1}^{i-1}\mathbbm{1}(B_j))\mathsf{E}_{\mathcal{F}_{t_{i-1}}}\mathbbm{1}(B_i^c)e^{\alpha (T_i - t_{i-1}+t_{i-1})}$
            
            $\stackrel{6}{\leq{}}
            \sum\limits_{i \geq{} 1}\mathsf{E}_{v_0}(\prod\limits_{j = 1}^{i-1}\mathbbm{1}(B_j))\mathsf{E}_{\mathcal{F}_{t_{i-1}}}\mathbbm{1}(B_i^c)e^{\alpha X_{t_{i-1}}}e^{\alpha t_{i-1}}$

            $\stackrel{7}{\leq{}}\sum\limits_{i \geq{} 1}\mathsf{E}_{v_0}(\prod\limits_{j = 1}^{i-1}\mathbbm{1}(B_j))(\mathsf{E}_{\mathcal{F}_{t_{i-1}}}\mathbbm{1}(B_i^c)e^{2\alpha X_{t_{i-1}}})^{1\over{2}}(\mathsf{E}_{\mathcal{F}_{t_{i-1}}}\mathbbm{1}(B_i^c)e^{2\alpha t_{i-1} })^{1\over{2}}$

            $ \stackrel{8}{\leq{}}\sum\limits_{i \geq{} 1}(\mathsf{E}_{v_0}\mathsf{E}_{\mathcal{F}_{t_{i-1}}}(\prod\limits_{j = 1}^{i-1}\mathbbm{1}(B_j))\mathbbm{1}(B_i^c)e^{2\alpha t_{i-1} })^{1\over{2}}(\mathsf{E}_{v_0}\mathsf{E}_{\mathcal{F}_{t_{i-1}}}(\prod\limits_{j = 1}^{i-1}\mathbbm{1}(B_j))\mathbbm{1}(B_i^c) e^{2\alpha X_{t_{i-1}}})^{1\over{2}}$
            
            $ = \sum\limits_{i \geq{} 1}(\mathsf{E}_{v_0}(\prod\limits_{j = 1}^{i-1}\mathbbm{1}(B_j))\mathbbm{1}(B_i^c)e^{2\alpha t_{i-1} })^{1\over{2}}(\mathsf{E}_{v_0}(\prod\limits_{j = 1}^{i-1}\mathbbm{1}(B_j))\mathbbm{1}(B_i^c) e^{2\alpha X_{t_{i-1}}})^{1\over{2}}.$
        \end{center}
        In (7) and (8) Cauchy-Bunyakovsky-Schwarz inequality was applied. We will evaluate separately $\mathsf{E}_{v_0}(\prod\limits_{j = 1}^{i-1}\mathbbm{1}(B_j))e^{2\alpha t_{i-1} }$ and $\mathsf{E}_{v_0}(\prod\limits_{j = 1}^{i-1}\mathbbm{1}(B_j))\mathbbm{1}(B_i^c) e^{2\alpha X_{t_{i-1}}}$.

        1) $$t_{i-1} = (t_{i-1} - T_{i-1}) + (T_{i-2} - t_{i-2})+\dots{}+(T_1-t_0)+(t_0-T_0).$$
        
        We have
        \begin{center}
        $ \mathsf{E}_{v_0}(\prod\limits_{j = 1}^{i-1}\mathbbm{1}(B_j))e^{2\alpha t_{i-1} } = 
        \mathsf{E}_{v_0}(\prod\limits_{j = 1}^{i-1}\mathbbm{1}(B_j))(\prod\limits_{k = 0}^{i-1}e^{2\alpha (t_{k} - T_k) })(\prod\limits_{l = 1}^{i-1}e^{2\alpha (T_{l} - t_{l-1}) })$
        $=\mathsf{E}_{v_0}\mathsf{E}_{\mathcal{F}_{T_{i-1}}}(\prod\limits_{j = 1}^{i-1}\mathbbm{1}(B_j))(\prod\limits_{k = 0}^{i-1}e^{2\alpha (t_{k} - T_k) })(\prod\limits_{l = 1}^{i-1}e^{2\alpha (T_{l} - t_{l-1}) })$
        $=\mathsf{E}_{v_0}(\prod\limits_{j = 1}^{i-1}\mathbbm{1}(B_j))(\prod\limits_{k = 0}^{i-2}e^{2\alpha (t_{k} - T_k) })(\prod\limits_{l = 1}^{i-1}e^{2\alpha (T_{l} - t_{l-1}) })\mathsf{E}_{\mathcal{F}_{T_{i-1}}}e^{2\alpha (t_{i-1} - T_{i-1}) } $
        $\stackrel{\text{lemma  1}}{\leq{}}M_1^{(2\alpha)}\mathsf{E}_{v_0}(\prod\limits_{j = 1}^{i-1}\mathbbm{1}(B_j))(\prod\limits_{k = 0}^{i-2}e^{2\alpha (t_{k} - T_k) })(\prod\limits_{l = 1}^{i-1}e^{2\alpha (T_{l} - t_{l-1}) }) $
        $=M_1^{(2\alpha)} \mathsf{E}_{v_0}\mathsf{E}_{\mathcal{F}_{t_{i-2}}}(\prod\limits_{j = 1}^{i-1}\mathbbm{1}(B_j))(\prod\limits_{k = 0}^{i-2}e^{2\alpha (t_{k} - T_k) })(\prod\limits_{l = 1}^{i-1}e^{2\alpha (T_{l} - t_{l-1}) })$
        $=M_1^{(2\alpha)} \mathsf{E}_{v_0}(\prod\limits_{j = 1}^{i-2}\mathbbm{1}(B_j))(\prod\limits_{l = 1}^{i-2}e^{2\alpha (T_{l} - t_{l-1}) })(\prod\limits_{k = 0}^{i-2}e^{2\alpha (t_{k} - T_k) }) \mathsf{E}_{\mathcal{F}_{t_{i-2}}} \mathbbm{1}(B_{i-1}) e^{2\alpha (T_{i-1} - t_{i-2}) }$
        $\stackrel{\text{lemma  2}}{\leq{}}
        M_1^{(2\alpha)} M_2^{(2\alpha)} \mathsf{E}_{v_0}(\prod\limits_{j = 1}^{i-2}\mathbbm{1}(B_j))(\prod\limits_{l = 1}^{i-2}e^{2\alpha (T_{l} - t_{l-1}) })(\prod\limits_{k = 0}^{i-2}e^{2\alpha (t_{k} - T_k) }). 
        $        
        \end{center}
        So, by induction we obtain
        \begin{center}
        $\mathsf{E}_{v_0}(\prod\limits_{j = 1}^{i-1}\mathbbm{1}(B_j))e^{2\alpha t_{i-1} }\leq{} (M_1^{(2\alpha)})^i (M_2^{(2\alpha)})^{i-1}.$
        \end{center}
        
        2) Note that $X_{t_i} \geq X_{T_i}$ и $X_{T_0} = X_{t_0}=x$, then
        \begin{center}
        $X_{t_{i-1}} = (X_{t_{i-1}}- X_{t_{i-2}}) + (X_{t_{i-2}}- X_{t_{i-3}})+\dots+(X_{t_{1}}- X_{t_{0}} ) +(X_{t_0} - x) + x$
        $\leq{}(X_{t_{i-1}}- X_{T_{i-1}}) + (X_{t_{i-2}}- X_{T_{i-2}})+\dots+(X_{t_{1}}- X_{T_{1}} ) +(X_{t_0} - X_{T_0}) + x.$
        \end{center}
        
        Let estimate
        \begin{center}
        $\mathsf{E}_{v_0}(\prod\limits_{j = 1}^{i-1}\mathbbm{1}(B_j))\mathbbm{1}(B_i^c) e^{2\alpha X_{t_{i-1}}} = \mathsf{E}_{v_0}(\prod\limits_{j = 1}^{i-1}\mathbbm{1}(B_j))\mathbbm{1}(B_i^c) e^{2\alpha x + 2\alpha(X_{t_0} - x) + 2\alpha \sum\limits_{k=1}^{i-1}(X_{t_k} - X_{t_{k-1}})}$
        $\leq e^{2\alpha x}\mathsf{E}_{v_0}(\prod\limits_{j = 1}^{i-1}\mathbbm{1}(B_j))e^{ 2\alpha \sum\limits_{k=0}^{i-1}(X_{t_k} - X_{T_k})} = e^{2\alpha x}\mathsf{E}_{v_0}\mathsf{E}_{\mathcal{F}_{T_{i-1}}}(\prod\limits_{j = 1}^{i-1}\mathbbm{1}(B_j))e^{ 2\alpha \sum\limits_{k=0}^{i-1}(X_{t_k} - X_{T_k})}$
        $= e^{2\alpha x}\mathsf{E}_{v_0}(\prod\limits_{j = 1}^{i-1}\mathbbm{1}(B_j))e^{ 2\alpha \sum\limits_{k=0}^{i-2}(X_{t_k} - X_{T_k})}\mathsf{E}_{\mathcal{F}_{T_{i-1}}}e^{ 2\alpha (X_{t_{i-1}} - X_{T_{i-1}})}$
        $\stackrel{\text{lemma  3}}{\leq} e^{2\alpha x} M_3^{(2\alpha)}\mathsf{E}_{v_0}(\prod\limits_{j = 1}^{i-1}\mathbbm{1}(B_j))e^{ 2\alpha \sum\limits_{k=0}^{i-2}(X_{t_k} - X_{T_k})}.$
        \end{center}

        By induction and considering that $M_3^{(2\alpha)} \ge 1$, we get

        \begin{center}
        $
        \mathsf{E}_{v_0}(\prod\limits_{j = 1}^{i-1}\mathbbm{1}(B_j))\mathbbm{1}(B_i^c) e^{2\alpha X_{t_{i-1}}}  \le e^{2\alpha x} (M_3^{(2\alpha)})^{i-1}\mathsf{E}_{v_0}(\prod\limits_{j = 1}^{i-1}\mathbbm{1}(B_j))e^{ 2\alpha (X_{t_0} - X_{T_0})} 
        $
        $=e^{2\alpha x} (M_3^{(2\alpha)})^{i-1}\mathsf{E}_{v_0}\mathsf{E}_{\mathcal{F}_{t_0}}(\prod\limits_{j = 1}^{i-1}\mathbbm{1}(B_j))e^{ 2\alpha (X_{t_0} - X_{T_0})}$
        $=e^{2\alpha x} (M_3^{(2\alpha)})^{i-1}\mathsf{E}_{v_0}e^{ 2\alpha (X_{t_0} - X_{T_0})}\mathsf{E}_{\mathcal{F}_{t_0}}(\prod\limits_{j = 1}^{i-1}\mathbbm{1}(B_j))\le e^{2\alpha x} (M_3^{(2\alpha)})^{i-1}\bar q^{i-1}\mathsf{E}_{v_0}e^{ 2\alpha (X_{t_0} - X_{T_0})}$
        $ = e^{2\alpha x} (M_3^{(2\alpha)})^{i-1}\bar q^{i-1}\mathsf{E}_{v_0}1
        \le e^{2\alpha x} (M_3^{(2\alpha)})^i \bar q^{i-1}.
        $
        \end{center}

       \textbf{Case \RNumb{2}}: at $t=0$ process is going up.

        Let us define stopping times:
        \begin{center}
            $T_0 = 0,$ $t_0 = \xi_{0},$ $ T_1 = \chi_{t_0},$ $t_1 = \xi_{T_1},$ $ T_2 = \chi_{t_1},$ $t_2 = \xi_{T_2},$ $ T_3 = \chi_{t_2}, \dots$
        \end{center}
        Like in case \RNumb{1}, $ T_i$ - the moment when the next fall of the process ends, which began in $t_{i-1}$, $ t_i$ accordingly, the moment when growth was replaced by a fall. Process a.s. finite number of times will
        change its behavior until it reaches $[0, N]$. Note that $T_i-t_{i-1}\leq{X_{t_{i-1}}}$ and $B_j \in \mathcal F_{T_j}$.

       As before
        \begin{center}
            $\mathsf{E}_{v_0}e^{\alpha \tau} =
            \sum\limits_{i \geq{} 1} \mathsf{E}_{v_0}e^{\alpha \tau}\mathbbm{1}(A_i)= 
            \sum\limits_{i \geq{} 1} \mathsf{E}_{v_0}e^{\alpha T_i}\mathbbm{1}((\cap_{j=1}^{i-1}B_j)\cap B_i^c)\stackrel{2}{=} 
            \sum\limits_{i \geq{} 1} \mathsf{E}_{v_0}(\prod\limits_{j=1}^{i-1}\mathbbm{1}(B_j))\mathbbm{1}(B_i^c)e^{\alpha T_i}$
            
            $\stackrel{3}{=}
            \sum\limits_{i \geq{} 1}\mathsf{E}_{v_0}\mathsf{E}_{\mathcal{F}_{t_{i-1}}}(\prod\limits_{j = 1}^{i-1}\mathbbm{1}(B_j))\mathbbm{1}(B_i^c)e^{\alpha T_i}\stackrel{4}{=}
            \sum\limits_{i \geq{} 1}\mathsf{E}_{v_0}(\prod\limits_{j = 1}^{i-1}\mathbbm{1}(B_j))\mathsf{E}_{\mathcal{F}_{t_{i-1}}}\mathbbm{1}(B_i^c)e^{\alpha T_i}$

            $\stackrel{5}{=}\sum\limits_{i \geq{} 1}\mathsf{E}_{v_0}(\prod\limits_{j = 1}^{i-1}\mathbbm{1}(B_j))\mathsf{E}_{\mathcal{F}_{t_{i-1}}}\mathbbm{1}(B_i^c)e^{\alpha (T_i - t_{i-1}+t_{i-1})}$
            
            $\stackrel{6}{\leq{}}
            \sum\limits_{i \geq{} 1}\mathsf{E}_{v_0}(\prod\limits_{j = 1}^{i-1}\mathbbm{1}(B_j))\mathsf{E}_{\mathcal{F}_{t_{i-1}}}\mathbbm{1}(B_i^c)e^{\alpha X_{t_{i-1}}}e^{\alpha t_{i-1}}$

            $\stackrel{7}{\leq{}}\sum\limits_{i \geq{} 1}\mathsf{E}_{v_0}(\prod\limits_{j = 1}^{i-1}\mathbbm{1}(B_j))(\mathsf{E}_{\mathcal{F}_{t_{i-1}}}\mathbbm{1}(B_i^c)e^{2\alpha X_{t_{i-1}}})^{1\over{2}}(\mathsf{E}_{\mathcal{F}_{t_{i-1}}}\mathbbm{1}(B_i^c)e^{2\alpha t_{i-1} })^{1\over{2}}$

            $ \stackrel{8}{\leq{}}\sum\limits_{i \geq{} 1}(\mathsf{E}_{v_0}\mathsf{E}_{\mathcal{F}_{t_{i-1}}}(\prod\limits_{j = 1}^{i-1}\mathbbm{1}(B_j))\mathbbm{1}(B_i^c)e^{2\alpha t_{i-1} })^{1\over{2}}(\mathsf{E}_{v_0}\mathsf{E}_{\mathcal{F}_{t_{i-1}}}(\prod\limits_{j = 1}^{i-1}\mathbbm{1}(B_j))\mathbbm{1}(B_i^c) e^{2\alpha X_{t_{i-1}}})^{1\over{2}}$
            
            $ = \sum\limits_{i \geq{} 1}(\mathsf{E}_{v_0}(\prod\limits_{j = 1}^{i-1}\mathbbm{1}(B_j))\mathbbm{1}(B_i^c)e^{2\alpha t_{i-1} })^{1\over{2}}(\mathsf{E}_{v_0}(\prod\limits_{j = 1}^{i-1}\mathbbm{1}(B_j))\mathbbm{1}(B_i^c) e^{2\alpha X_{t_{i-1}}})^{1\over{2}}.$
        \end{center}
        In (7) and (8) Cauchy-Bunyakovsky-Schwarz inequality was applied. We will evaluate separately $\mathsf{E}_{v_0}(\prod\limits_{j = 1}^{i-1}\mathbbm{1}(B_j))e^{2\alpha t_{i-1} }$ и $\mathsf{E}_{v_0}(\prod\limits_{j = 1}^{i-1}\mathbbm{1}(B_j))\mathbbm{1}(B_i^c) e^{2\alpha X_{t_{i-1}}}$.

        1) $$t_{i-1} = (t_{i-1} - T_{i-1}) + (T_{i-2} - t_{i-2})+\dots{}+(T_1-t_0)+(t_0-T_0)$$

        \begin{center}
        $ \mathsf{E}_{v_0}(\prod\limits_{j = 1}^{i-1}\mathbbm{1}(B_j))e^{2\alpha t_{i-1} } = 
        \mathsf{E}_{v_0}(\prod\limits_{j = 1}^{i-1}\mathbbm{1}(B_j))(\prod\limits_{k = 0}^{i-1}e^{2\alpha (t_{k} - T_k) })(\prod\limits_{l = 1}^{i-1}e^{2\alpha (T_{l} - t_{l-1}) }) =$
        $\mathsf{E}_{v_0}\mathsf{E}_{\mathcal{F}_{T_{i-1}}}(\prod\limits_{j = 1}^{i-1}\mathbbm{1}(B_j))(\prod\limits_{k = 0}^{i-1}e^{2\alpha (t_{k} - T_k) })(\prod\limits_{l = 1}^{i-1}e^{2\alpha (T_{l} - t_{l-1}) })$
        $=\mathsf{E}_{v_0}(\prod\limits_{j = 1}^{i-1}\mathbbm{1}(B_j))(\prod\limits_{k = 0}^{i-2}e^{2\alpha (t_{k} - T_k) })(\prod\limits_{l = 1}^{i-1}e^{2\alpha (T_{l} - t_{l-1}) })\mathsf{E}_{\mathcal{F}_{T_{i-1}}}e^{2\alpha (t_{i-1} - T_{i-1}) } $
        $\stackrel{\text{lemma 1}}{\leq{}}M_1^{(2\alpha)}\mathsf{E}_{v_0}(\prod\limits_{j = 1}^{i-1}\mathbbm{1}(B_j))(\prod\limits_{k = 0}^{i-2}e^{2\alpha (t_{k} - T_k) })(\prod\limits_{l = 1}^{i-1}e^{2\alpha (T_{l} - t_{l-1}) })$
        $=M_1^{(2\alpha)} \mathsf{E}_{v_0}\mathsf{E}_{\mathcal{F}_{t_{i-2}}}(\prod\limits_{j = 1}^{i-1}\mathbbm{1}(B_j))(\prod\limits_{k = 0}^{i-2}e^{2\alpha (t_{k} - T_k) })(\prod\limits_{l = 1}^{i-1}e^{2\alpha (T_{l} - t_{l-1}) })$
        $=M_1^{(2\alpha)} \mathsf{E}_{v_0}(\prod\limits_{j = 1}^{i-2}\mathbbm{1}(B_j))(\prod\limits_{l = 1}^{i-2}e^{2\alpha (T_{l} - t_{l-1}) })(\prod\limits_{k = 0}^{i-2}e^{2\alpha (t_{k} - T_k) }) \mathsf{E}_{\mathcal{F}_{t_{i-2}}} \mathbbm{1}(B_{i-1}) e^{2\alpha (T_{i-1} - t_{i-2}) }$
        $\stackrel{\text{lemma 2}}{\leq{}}
        M_1^{(2\alpha)} M_2^{(2\alpha)} \mathsf{E}_{v_0}(\prod\limits_{j = 1}^{i-2}\mathbbm{1}(B_j))(\prod\limits_{l = 1}^{i-2}e^{2\alpha (T_{l} - t_{l-1}) })(\prod\limits_{k = 0}^{i-2}e^{2\alpha (t_{k} - T_k) }).
        $        
        \end{center}
        So, by induction we obtain 
        \begin{center}
        $\mathsf{E}_{v_0}(\prod\limits_{j = 1}^{i-1}\mathbbm{1}(B_j))e^{2\alpha t_{i-1} }\leq{} (M_1^{(2\alpha)})^i (M_2^{(2\alpha)})^{i-1}.$
        \end{center}
        
        2) Note that $X_{t_i} \geq X_{T_i}$ and $X_{T_0} = X_{t_0}=x$, then
        \begin{center}
        $X_{t_{i-1}} = (X_{t_{i-1}}- X_{t_{i-2}}) + (X_{t_{i-2}}- X_{t_{i-3}})+\dots+(X_{t_{1}}- X_{t_{0}} ) +(X_{t_0} - x) + x$
        $\leq{}(X_{t_{i-1}}- X_{T_{i-1}}) + (X_{t_{i-2}}- X_{T_{i-2}})+\dots+(X_{t_{1}}- X_{T_{1}} ) +(X_{t_0} - X_{T_0}) + x$
        \end{center}
        
        further,
        \begin{center}
        $\mathsf{E}_{v_0}(\prod\limits_{j = 1}^{i-1}\mathbbm{1}(B_j))\mathbbm{1}(B_i^c) e^{2\alpha X_{t_{i-1}}} = \mathsf{E}_{v_0}(\prod\limits_{j = 1}^{i-1}\mathbbm{1}(B_j))\mathbbm{1}(B_i^c) e^{2\alpha x + 2\alpha(X_{t_0} - x) + 2\alpha \sum\limits_{k=1}^{i-1}(X_{t_k} - X_{t_{k-1}})}$
        $\leq e^{2\alpha x}\mathsf{E}_{v_0}(\prod\limits_{j = 1}^{i-1}\mathbbm{1}(B_j))e^{ 2\alpha \sum\limits_{k=0}^{i-1}(X_{t_k} - X_{T_k})} = e^{2\alpha x}\mathsf{E}_{v_0}\mathsf{E}_{\mathcal{F}_{T_{i-1}}}(\prod\limits_{j = 1}^{i-1}\mathbbm{1}(B_j))e^{ 2\alpha \sum\limits_{k=0}^{i-1}(X_{t_k} - X_{T_k})}$
        $= e^{2\alpha x}\mathsf{E}_{v_0}(\prod\limits_{j = 1}^{i-1}\mathbbm{1}(B_j))e^{ 2\alpha \sum\limits_{k=0}^{i-2}(X_{t_k} - X_{T_k})}\mathsf{E}_{\mathcal{F}_{T_{i-1}}}e^{ 2\alpha (X_{t_{i-1}} - X_{T_{i-1}})}$
        $\stackrel{\text{lemma 3}}{\leq} e^{2\alpha x} M_3^{(2\alpha)}\mathsf{E}_{v_0}(\prod\limits_{j = 1}^{i-1}\mathbbm{1}(B_j))e^{ 2\alpha \sum\limits_{k=0}^{i-2}(X_{t_k} - X_{T_k})}.$
        \end{center}

        By induction we get

        \begin{center}
        $
        \mathsf{E}_{v_0}(\prod\limits_{j = 1}^{i-1}\mathbbm{1}(B_j))\mathbbm{1}(B_i^c) e^{2\alpha X_{t_{i-1}}}  \le e^{2\alpha x} (M_3^{(2\alpha)})^{i-1}\mathsf{E}_{v_0}(\prod\limits_{j = 1}^{i-1}\mathbbm{1}(B_j))e^{ 2\alpha (X_{t_0} - X_{T_0})} 
        $
        $=e^{2\alpha x} (M_3^{(2\alpha)})^{i-1}\mathsf{E}_{v_0}\mathsf{E}_{\mathcal{F}_{t_0}}(\prod\limits_{j = 1}^{i-1}\mathbbm{1}(B_j))e^{ 2\alpha (X_{t_0} - X_{T_0})}$
        $=e^{2\alpha x} (M_3^{(2\alpha)})^{i-1}\mathsf{E}_{v_0}e^{ 2\alpha (X_{t_0} - X_{T_0})}\mathsf{E}_{\mathcal{F}_{t_0}}(\prod\limits_{j = 1}^{i-1}\mathbbm{1}(B_j))\le e^{2\alpha x} (M_3^{(2\alpha)})^{i-1}\bar q^{i-1}\mathsf{E}_{v_0}e^{ 2\alpha (X_{t_0} - X_{T_0})}
        \stackrel{\text{lemma 3}}{\le} e^{2\alpha x} (M_3^{(2\alpha)})^i \bar q^{i-1}.
        $
        \end{center}

        In both cases we get the same estimation
        \begin{center}
        $\mathsf{E}_{v_0}e^{\alpha \tau} \le \sum\limits_{i \geq{} 1}(\mathsf{E}_{v_0}(\prod\limits_{j = 1}^{i-1}\mathbbm{1}(B_j))\mathbbm{1}(B_i^c)e^{2\alpha t_{i-1} })^{1\over{2}}(\mathsf{E}_{v_0}(\prod\limits_{j = 1}^{i-1}\mathbbm{1}(B_j))\mathbbm{1}(B_i^c) e^{2\alpha X_{t_{i-1}}})^{1\over{2}}$
        $\le \sum\limits_{i \geq{} 1}((M_1^{(2\alpha)})^i (M_2^{(2\alpha)})^{i-1})^{1\over{2}} (e^{2\alpha x} (M_3^{(2\alpha)})^i \bar q^{i-1})^{1\over{2}}= {e^{\alpha x} \sqrt{M_1^{(2\alpha)}M_3^{(2\alpha)}}\over{1- \sqrt{M_1^{(2\alpha)}M_2^{(2\alpha)}M_3^{(2\alpha)}\bar q}}}.$
        \end{center}
         Theorem is proved. \hfill $\square$ 

    \section{Acknowledgments}
    The author is grateful to the supervisor Professor Alexander Veretennikov for the problem and the helpful discussions concerning this paper.

\end{document}